\documentclass[a4paper]{paper}
\usepackage{latexsym}
\usepackage{amsthm,amsmath,amssymb}

\usepackage[dvips]{graphicx, color}
\usepackage[title]{appendix}
 \usepackage{color}

\usepackage{color}
\usepackage{xcolor}
\usepackage{bbm}
\usepackage{scalerel}
\usepackage{stackengine}

\newtheorem{thm}  {Theorem}

\newtheorem*{theorem-non}  {Theorem}

\newtheorem{prop}[thm]{Proposition}

\newtheorem*{condition-non}  {Condition}
\newtheorem{ex}[]{Example}

\newtheorem{cor}[thm]{Corollary}

\newtheorem{lem}[thm]{Lemma}

\newtheorem{defn}[thm]{Definition}
\newtheorem{rem}[]{Remark}

\newcommand{\cE}  {{\mathcal E}}

\newcommand{\cG}  {{\mathcal G}}

\newcommand{\cL}  {{\mathcal L}}
\newcommand{\cM}  {{\mathcal M}}

\newcommand{\A}  {{\Bbb A}}

\newcommand{\C}  {{\Bbb C}}

\newcommand{\F}  {{\Bbb F}}
\newcommand{\G}  {{\Bbb G}}

\renewcommand{\O}  {{\Bbb O}}

\newcommand{\Q}  {{\Bbb Q}\hspace{.06em}}
\newcommand{\R}  {{\Bbb R}}

\newcommand{\Z}  {{\Bbb Z}}

\newcommand{\frb}  {{\frak  b}}
\newcommand{\frc}  {{\frak  c}}

\newcommand{\frg}  {{\frak  g}}

\newcommand{\frk}  {{\frak  k}}
\newcommand{\frl}  {{\frak  l}}

\newcommand{\frp}  {{\frak  p}}

\newcommand{\frs}  {{\frak  s}}

\def\O{\mathcal O}

\def\ka{\kappa}
\def\lla{\longlefttarrow}
\def\lra{\longrightarrow}

\def\={\:=\:}  \def\+{\,+\,}

\def\a{\alpha} \def\b{{\beta}}  \def\ba{\overline\a}

\newcommand{\frt}  {{\frak  t}}

\newcommand{\frv}  {{\frak  v}}

\newcommand{\frz}  {{\frak  z}}

\newcommand{\frP}  {{\frak  P}}

\newcommand{\frS}  {{\frak  S}}

\def\be{\begin{equation}}   \def\ee{\end{equation}}
\def\bes{\begin{equation*}}   \def\ees{\end{equation*}}
\def\ba{\begin{aligned}}   \def\ea{\end{aligned}}
\def\bc{\begin{cases}}   \def\ec{\end{cases}}
\def\bp{\begin{proof}}   \def\ep{\end{proof}}

\newcommand{\ar}  {\mathrm{ar}}

\newcommand{\fin}  {\mathrm{fin}}

\newcommand{\Aut}  {\mathrm{Aut}}

\newcommand{\sep}  {\mathrm{\sep}}

\newcommand{\ov}  {\overline}
\def\qqan{\qquad\mathrm{and}\qquad}
\def\qan{\quad\mathrm{and}\quad}

\newcommand{\dis}  {\displaystyle}

\def\smm{\smallsetminus}
\def\s{\sigma}
\def\La{\Lambda}
\def\la{\lambda}

\def\Ker{\mathrm{Ker}}

\def\GL{\mathrm{GL}}

\def\lra{\longrightarrow}
\def\lla{\longlefttarrow}

\def\ov{\overline}
\def\De{\Delta}
\def\lan{\langle}
\def\ran{\rangle}

\def\lla{\longleftarrow}

\def\lan{\langle}
\def\ran{\rangle}

\def\bbm1{\mathbbm_1}

\def\ka{\kappa}
\def\De{\Delta}

\def\be{\begin{equation}}   \def\ee{\end{equation}}
\def\bes{\begin{equation*}}   \def\ees{\end{equation*}}
\def\bea{\begin{equation}\begin{aligned}}   
\def\eea{\end{aligned}\end{equation}}

\def\tot{\mathrm{tot}}
\def\Lie{\mathrm{Lie}}
\def\Ad{\mathrm{Ad}}
\def\ad{\mathrm{ad}}

\def\bm{\begin{matrix}}
\def\em{\end{matrix}}
\def\bpm{\begin{pmatrix}}
\def\epm{\end{pmatrix}}
\def\Hom{\mathrm{Hom}}

\def\rk{\mathrm{rank}}
\def\diag{\mathrm{diag}}

\def\Spec{\mathrm{Spec}}

\def\Gal{\mathrm{Gal}}

\def\End{\mathrm{End}}

\def\Out{\mathrm{Out}}
\def\h2m{\hskip 2.0cm}

\def\ed{\end{document}}

\def\sep{\mathrm{sep}}

\def\ses{\mathrm{ss}}
\def\ab{\mathrm{ab}}

\newcommand\xdownarrow[1][2ex]{%
   \mathrel{\rotatebox{90}{$\xleftarrow{\rule{#1}{0pt}}$}}
}

\begin{document}
\title{\bf Arithmetic Characteristic Curves} 
\author{Lin WENG}  
\date{}
\maketitle
\begin{abstract}
For a split reductive group defined over a number field, we first 
introduce the notations of arithmetic torsors and arithmetic Higgs 
torsors. Then we construct  arithmetic characteristic curves associated 
to arithmetic Higgs torsors, based on the Chevalley characteristic 
morphism and the existence of Chevalley basis for the associated Lie 
algebra. As to be expected, this work is motivated by the works of 
Beauville-Narasimhan on spectral curves and Donagi-Gaistgory on 
cameral curves in algebraic geometry. In the forthcoming papers, we 
will use arithmetic characteristic curves to construct arithmetic Hitchin 
fibrations and study the intersection homologies and 
perverse sheaves for the associated structures, following Ngo's 
approach to the fundamental lemma.
\end{abstract}

\section{Chevelley's Characteristic Morphism}
\subsection{Over Number Fields}
Let $F$ be a number field with $\O_F$ its ring of integers. Denote by 
$X=\Spec\,\O_F$ the associated uncompleted arithmetic curve. 

Let $G$ be a split reductive group over $F$. Fix a split maximal 
subtorus $T$ and a maximal split quotient torus $T'$ of $G$. Denote 
the Lie algebra of $G$ by $\frg:=\Lie\,G$ , and set $\frt:=\Lie\, T$ be 
the associated commutative subalgebra of $\frg$.

Recall that, with respect to the adjoint action 
\bes
\bm
\ad: &\frg\times\frg&\lra& \frg\\[0.60em]
&(g,x)&\mapsto&(\ad(g))(x):=[g,x]
\em
\ees
$\frg$ admits a natural decomposition
\bes
\frg:=\frt\,\bigoplus \bigoplus_\a\frg_\a
\ees
where 
\bes
\frg_\a:=\big\{x\in\frg: (\ad(h))(x)=\a(h)\,x\big\}
\ees
for $\a$ running through a finite subset $\Phi$ of the space 
\bes
X^*(T):=\Hom(T,\G_m),
\ees
of rational characteristics of $T$.

For a fixed minimal split parabolic subgroup $B$ of $G$ containing $T$, 
set $\frb:=\Lie\,B$. Then $G/B$ is proper, and there exists a finite 
subset $\Phi_+$ of $\Phi$, the so-called set of positive roots of 
$(G,B,T)$,  such that
\begin{enumerate}
\item[(1)] $\Phi=\Phi_+\bigsqcup\, (-\Phi_+)$, 
\item[(2)] $\frt\oplus\oplus^{~}_{\a\in \Phi_+}\frg_\a\subset \frb$, and 
\item[(3)] $\Phi_+$ admits a subset $\De$ of simple roots associated 
to $(G,B,T)$, such that
\begin{enumerate}
\item [(i)] $\Phi_+\subseteq\sum_{\a\in\De}\Z_{\geq 0}\a$, where $\Z_{\geq 0}:=\big\{n\in\Z:\,n\geq 0\big\}$, and
\item[(ii)] $\De$ forms a basis of the $\Q$-linear space 
\bes
X^*(T)^{~}_\Q:=X^*(T)\otimes\Q.
\ees  
\end{enumerate}
\end{enumerate}

Let $W$ be the Weyl group of $G$, defined as the finite quotient group 
\bes
W:=N_G(T)/Z_G(T)
\ees
where $N_G(T)$, resp. $Z_G(T)$, denotes the normalizer subgroup, resp. the centralizer subgroup, of $T$ in $G$. It is well known that $W$ is canonically isomorphic to the subgroup of the automorphism group of $X^*(T)_\Q$ generated by the reflections
\bes
\bm
\s_\a:&X^*(T)^{~}_\Q&\lra& X^*(T)^{~}_\Q\\[0.5em]
&v&\mapsto& v-\frac{2}{(v,a)}\,\a
\em
\ees
It is a canonical result due to Chevalley that, over the base field $F$, the 
space of $G$-invariant polynomials on $G$ coincides with the space of 
the $W$-invariant polynomials of $T$. That is to say,
\be\label{eq1}
F[G]^G\simeq F[T]^W
\ee
where the actions of both sides are defined by  
\bes
\bm G\times F[G]&\lra& F[G]\\[0.6em]
(g,\sum_{i} a_ig_i)&\mapsto&\sum_ia_i(gg_ig^{-1})\em\qan
\bm W\times F[T]&\lra& F[T]\\[0.6em]
(\s,\sum_{i} a_it_i)&\mapsto&\sum_ia_i\s(t_i)
\em
\ees
Similarly, in terms of the Lie algebras, we have
\bes
F[\frg]^G\simeq F[\frt]^W
\ees
where $G$ acts on $\frg$ in terms of the adjoint action $\Ad$, namely, 
$\Ad(g)$ is defied as the differential of the conjugation morphism 
$x\mapsto gxg^{-1}$ of $G$. 

In terms of geometry, the isomorphism \eqref{eq1} naturally induces 
the  scheme theoretic morphism
\bes
\bm
&&\Spec\,F[T]\\
&&\big\downarrow\\
G&\buildrel{\chi}\over\lra&T/\!/W:=\Spec\,F[T]^W
\em\ees
or equivalently, for the associated Lie structures,
\be\label{eq2} 
\bm
\frg_F\simeq\Spec\,F[\frg]& &\frt_F\simeq\Spec\,F[\frt]\\
\Big\downarrow &\chi\searrow{~}&\Big\downarrow\pi_F^{~} \\
\frg_F/\!/G:=\Spec\,F[\frg]^G&\simeq&\frt_F/\!/W:=\Spec\,F[\frt]^W
\em 
\ee
where $\frg_F:=\frg\otimes F$ and similarly $\frt_F:=\frt\otimes F$. 

\begin{ex} \normalfont For $G=\GL_n/F$, we have $W\simeq\frS_n$, 
the symmetric group on $n$ symbols and 
$\frg_F=\frg\frl_n(F)=\End(F^n)$. Then $\chi$ coincides with the 
morphism
\bes
\bm
\chi:&\End(F^n)&\lra&\bigoplus_{k=1}^nL_k\\[0.5em]
&A&\mapsto&\det(\la I_n-A)
\em
\ees
where $L_k$ denotes the one dimensional vector space generated by 
the $k$-th elementary symmetric polynomials, and $I_n$ denotes the 
unity matrix of size $n$. That is to say, $\chi$ assigns a matrix $A$ to 
the associated eigen polynomial. In particular, by restricting $\chi$ to 
$\frt_F$ which consists of diagonal matrices $D=\diag(a_1,\ldots,a_n)$, 
we conclude that 
\bes
\chi(D)
=t^n-\sum_{i=1}^na_it^{n-1}+\sum_{i<j} a_i a_j t^2+\ldots+
(-1)^n\prod_{i=1}^na_i
\ees
For this reason, the morphism $\chi$ for general $G$ is called 
{\it the Chevelley characteristic morphism}. From above, after 
restricting $\chi$ to $\frt$, the Chevelley characteristic morphism 
is simply equivalent to the assignments of the unordered  eigenvalues.
\end{ex} 

\subsection{Over Integral Bases}

The diagram in \eqref{eq2} associated to a split reductive group 
$G/F$ only works over the pointed base $\Spec \,F$.
In this section, we construct a natural extension to the integral base 
$\Spec\,\O_F$, whose generic fiber (over the generic 
point $\eta_F:=\Spec\,F$) coincides with that of \eqref{eq2}.\\

To start with, we recall the so-called Chevelley basis for $\frg_F$.
For simplicity, we assume that $F=\Q$ for the time being.
 
\begin{defn}[See e.g. Ch. VII, \S25 of \cite{H}]\label{def1} \normalfont 
A {\it Chevelley basis} for $[\frg,\frg]$ is a basis for the 
$\Q$-linear space $[\frg,\frg]$, consisting of 
$\big\{x_\a:\a\in\Phi\big\}\bigsqcup\big\{h_\a:\a\in\De\big\}$ which 
satisfy the following properties:
\begin{enumerate}
\item [(a)] for all $\a\in \Phi$, $x_\a\in\frg_\a$.
\item[(b)]  for all $\a\in \Phi$, $[x_\a,x_{-\a}]=h_\a$ so that 
$x_\a, x_{-\a}$ and $h_\a$ span a three dimensional simple subalgebra 
of $\frg$ which is isomorphic to $\frs\frl_2(F)$ via
\bes
x_\a\mapsto \bpm 0&1\\0&0\epm,\quad x_{-\a}\mapsto \bpm 0&0\\1&0\epm, h_\a\mapsto \bpm 1&0\\0&1\epm
\ees
\item[(c)] for $\a,\,\b\in \Phi$, if $[x_\a,x_\b]=c_{\a,\b}x_{\a+\b}$, then 
\begin{enumerate}
\item [(i)] $c_{\a,\b}=c_{-\a,-\b}$.
\item [(ii)] $\dis{c_{\a,\b}^2=\ka\,(\ell+1)\frac{\,(\a+\b,\a+\b)\,}{(\b,\b)}}$ where the constants $\ell,\,\kappa$ are defined by the $\a$-string 
$\b-\ell\a,\cdots,\b+\ka\a$ through $\b$.
\end{enumerate}
\end{enumerate}
\end{defn}

We have the following well-known

\begin{thm}[Chevelley] Over $[\frg,\frg]$, we have
\begin{enumerate}
\item [(1)] There always exists a Chevelley basis 
$\big\{x_\a:\a\in\Phi\big\}\bigsqcup\big\{h_\a:\a\in\De\big\}$ on 
$[\frg,\frg]$,
\item[(2)] All the structural constants  lie in $\Z$. That is to say,
\begin{enumerate}
\item [(i)] for all $\a,\,\b\in\De$, $[h_\a,h_\b]=0$.
\item [(ii)]  for all $\a\in\Phi,\ \b\in\De$, 
$[h_\b,x_\a]=\lan\a,\b\ran x_\a$, where
$\dis{\lan \a,\b\ran:=2\frac{\,(\a,\b)\,}{(\b,\b)}}$.
\item [(iii)] for $\a\in\Phi$, $[x_\a, x_{-\a}]$ is a $\Z$-linear combination 
of $h_\a$\!\!'s  ($\a\in \De$).
\item [(iv)] If $\a, \b$ are independent roots and 
$\b-\ell\a,\cdots,\b+\ka\a$ is the $\a$-string through $\b$, then 
$[x_\a,x_\b]=\bc \qquad\  0&\ka=0,\\
\pm(\ell+1)\,x_{\a+\b}&\a+\b\in\Phi.\ec$
\end{enumerate}
\end{enumerate}
\end{thm}

Obviously, once $\De$ is fixed, $h_\a$ are uniquely determined if 
$\a\in\De$. In addition, for a general $\a\in\Phi$, if $x_\a$ is 
replaced by $c_\a x_a$, deduced from the conditions in (c) of 
Definition\,\ref{def1}, $\{c_\a\}_{\a\in \Phi}$ are bounded by the 
constrains:
\begin{enumerate}
\item [(i)] for all $\a\in \Phi$, $c_\a c_{-\a}=1$,
\item[(ii)] for all $\a,\,\b\in \Phi$, if $\a+\b\in\Phi$, then
$c_\a c_\b=\pm c_{\a+\b}$.
\end{enumerate}

Conversely, it is clear that if $\{c_\a\}_{\a\in \Phi}$ satisfies (these two 
conditions)i) and (ii) just mentioned, then 
$\big\{x_\a:\a\in\Phi\big\}\bigsqcup\big\{h_\a:\a\in\De\big\}$ forms a 
Chevelley basis of $[\frg,\frg]$ as well.\\

To treat the Lie algebra $\frg$ associated to the split reductive group $G/\Q$, it suffices to use the decomposition
\bes
\frg=[\frg,\frg]\oplus\frz
\ees
where $\frz$ denotes the center of $\frg$. Obviously, the integral 
bases for $\frz$ are parametrized by $\GL_{\dim_\Q\frz}(\Z)$. 
Hence, it is natural to define {\it a Chevelley basis of $\frg$} to be
the union of a Chevelley basis of $[\frg,\frg]$ as above and an integral 
basis of $\frz$. Denote by $\frg_\Z^{~}$, resp. $\frt_\Z^{~}$, resp. 
$\frz_\Z^{~}$, the associated lattice of $\frg$, resp. of $\frt$, resp. 
of $\frz$. Obviously, $\frg_\Z^{~}$ admits a natural Lie structure and 
does not depend on the chosen integral Chevalley basis.

Moreover, working with $\frg_\Z$, we obtain the following structural 
diagram over $\Z$:
\bes 
\bm
\frg_\Z^{~}\simeq\Spec\,\Z[\frg_\Z^{~}]& &\frt_\Z^{~}\simeq\Spec\,
\Z[\frt_\Z^{~}]\\
\Big\downarrow &\chi_\Z^{~}\searrow{~}&\Big\downarrow\pi_\Z^{~} \\
\frg_\Z^{~}/\!/G(\Z):=\Spec\,\Z[\frg_\Z^{~}]^{G(\Z)}
&\simeq&\frt_\Z^{~}/\!/W:=\Spec\,\Z[\frt_\Z^{~}]^W.
\em 
\ees
Obviously,  associated to the base change 
$\Spec\,F\hookrightarrow\Spec\,\Z$, we recover the diagram in 
\eqref{eq2}.

The same construction works  for a general number field $F$ instead of 
$\Q$. That is to say, we may use the base change 
$\Spec\,\O_F\lra\Spec\,Z$ to obtain the following diagram over the 
integral base $\Spec\,\O_F$:
\be\label{eq3} 
\bm
\frg_{\O_F}^{~}\simeq\Spec\,\O_F[\frg_{\O_F}^{~}]& &\frt_{\O_F}^{~}
\simeq\Spec\,\Z[\frt_\Z^{~}]\\
\Big\downarrow &\!\!\chi_{\O_F}^{~}\searrow\!\!
&\Big\downarrow\pi_{\O_F}^{~} \\
\frg_{\O_F}^{~}/\!/G({\O_F}):=\Spec\,{\O_F}[\frg_{\O_F}^{~}]^{G({\O_F})}
&\simeq&\frt_{\O_F}^{~}/\!/W:=\Spec\,{\O_F}[\frt_{\O_F}^{~}]^W.
\em 
\ee
whose general fiber is the similar diagram over $\Spec\,F$, which can 
be obtained from \eqref{eq3} by replacing the integer ring $\O_F$ with 
its associated number field $F$. 

To end this discussion, via the Minkowski embedding 
$\O_F\hookrightarrow F\hookrightarrow F_\infty:=\R^{r_1}\times 
\C^{r_2}$, we obtain the Lie algebra $\frg_\infty:=\frg\otimes\F_\infty$ 
and similarly $\frt_\infty$ and $\frz_\infty$. Furthermore, induced from 
the base change $\Spec\,F_\infty\to \Spec\,\O_F$, we obtain the 
following diagram on $F_\infty$:
\be\label{eq4} 
\bm
\frg_\infty^{~}\simeq\Spec\,F_\infty[\frg_\infty^{~}]& &\frt_\infty^{~}
\simeq\Spec\,F_\infty[\frt_\infty^{~}]\\
\Big\downarrow &\chi_\infty^{~}\searrow{~}
&\Big\downarrow\pi_\infty^{~} \\
\frg_\infty^{~}/\!/
G(F_\infty):=\Spec\,F_\infty[\frg_\infty^{~}]^{G(F_\infty)}
&\simeq&\frt_\infty^{~}/\!/W:=\Spec\,F_\infty[\frt_\infty^{~}]^W.
\em 
\ee
In particular, similar to the morphisms $\pi_F^{~}$ and $\pi_{\O_F}^{~}$, 
$\pi_\infty$ is a finite $(|W|:1)$-morphism, even supposed to be highly 
ramified in general. For later use, we denote 
$\frt_{\O_F}^{~}/\!/W,\,\frt_F^{~}/\!/W,$ and $\frt_\infty^{~}/\!/W$ by 
$\frc_{\O_F}^{~}, \,\frc_F^{~}$ and $\frc_\infty^{~}$, respectively.

\section{$G$-Torsors on $\ov{\Spec\,\O_F}$}

Let $F$ be a number field with $\O_F$ its ring of integers and $\A_F$ 
its adelic ring. Denote by $S$ the set of inequivalent normalized 
valuations of $F$, and by $S_\fin$, resp. $S_\infty$, the subsets of $S$ 
consisting of non-archimedean, resp. archimedean, valuations.
For each $v\in S$, denote by $F_v$ the $v$-completion of $F$. When 
$v\in S_\infty$, $F_v$ is isomorphic to either $\R$ or $\C$; and if 
$v\in S_\fin$, $F_v$ is a discrete valuation fields. Accordingly, the 
valuation is called real, or complex, or $v$-adic. For $v$-adic valuations, 
denote by $\O_v$  the associated valuation ring of $F_v$, by $\frP_v$ 
its maximal ideal, and by $k_v:=\O_v/\frP_v$ its residue field. It is 
well-known that $\O_v$, being a discrete valuation ring, admits only two 
prime ideals, namely, $\frP_v$ and $\{0\}$, and $k_v$ is a finite 
extension of $\F_p$, for a certain prime number $p\in \Z$. We call 
$[k_v:\F_p]$ the residue extension degree of $v$.  Based on all these, 
we have $\A_F=\prod_{v\in S}'F_v$ where $\prod_{v\in S}'$ denotes the 
restricted product of the $F_v$'s with respect to the $\O_c$'s. That is 
to say, an element $a=(a_v)\in\prod_{v\in S}F_v$ belongs to $\A_F$ if 
and only if $a_v\in \O_v$ for all but finitely many $v\in S$. It is well 
known that, induced from the locally compact topologies on the 
$F_v$'s, $\A_F$ is locally compact. 

Let $G$ be a split reductive group over $F$ with a pinning 
$(T,B,\{x_\a\}_{\a\in\De})$. Here $T$ is a maximal split subtorus of $G$, 
$B$ is a minimal split parabolic subgroup of $G$ containing $T$,  
and $\De$ denotes the set of simple roots of the root system
$\Phi$ associated to $(G,T,B)$ and  $x_\a$ denotes a non-zero vector 
of the proper subspace $\Lie(U)_\a$ of the Lie algebra $\Lie(U)$  
corresponding to the eigenvalue $\a$, where $U$ denotes the 
unipotent radical of $B$. We assume that $\{x_\a\}_{\a\in\De}$ can be 
extended to a Chevelley basis of $\frg^{~}_F$. Set then
\be
x_+:=\sum_{\a\in\De}x_\a.
\ee

\subsection{Torsors over Local and Global Fields}

Let $\ov F$ be the algebraic closure on $F$ contained in $\C$, and set 
$G_F=\Gal(\ov F/F)$ be the absolute Galois group of $F$.  For each 
$v\in S$, fix an algebraic closure $\ov F_v$ of $F_v$, and denote by 
$G_{F_v}=\Gal(\ov F_v/F_v)$ the local absolute Galois group of $F$ at 
$v$.

\begin{defn} \normalfont
For $K=F$ or $F_v$, a $K$-scheme $\cG$ is called a 
{\it $G$-torsor} if $\cG$ is equipped with a faithful, transitive and 
$G_K$-compatible action of $G(\ov K)$ on $\cG(\ov K)$.
\end{defn}

\begin{ex} For $K=F$ or $F_v$, set $\cG=G$. Then 
\bes
\bm
G(\ov K)\times \cG(\ov K)\times \Gal(\ov K/K)&\lra&\cG(\ov K)\\[0.6em]
(g,d,\s)&\mapsto&(g\,d)^\s=g^\s \,d^\s
\em
\ees
gives a $G$-torsor structure on $\cG$ over $K$.
\end{ex}

Since $G$ is defined over $F$, there is a continuous action of 
$G_F$ on $G(\ov F)$. As usual, set $H^0(F,G):=G(\ov F)^{\Gal(\ov F/F)}$ 
be the collection of $G_F$-invariant points of $G(\ov F)$ and 
$H^1(F,G):=Z^1(F,G)/B^1(F,G)$ be the set of equivalence classes of 
1-cocycles, where $Z^1(F,G)$, resp. $B^1(F,G)$, denotes the set of 
1-cocycles, resp. 1-coboundaries, of $G_F$ on $G(\ov F)$, i.e. a 
continuous map 
\bes
\bm
\phi:&G_F&\lra&G(\ov F)\\[0.6em]
&\s&\mapsto &a_\s
\em
\ees
satisfying
\bes
a_{\a\tau}=a_\s\cdot a_\tau^\s,
\ees
and two 1-cocycles $a_\s$ and $a_s'$ are said to be equivalence if 
there exists an element $g$ of 
$G(\ov F)$ such that
\bes
a_\s'=g^{-1}\cdot a_\s\cdot g^\s
\ees
In other words, $a_\s$ is an 1-coboundary if $\a_\s=g^{-1}g^\s$.

\begin{thm}[See e.g. \S2 of \cite{NSW}] For $K=F$ or $F_v$, there 
exists a natural bijection between the set of isomorphism classes of 
$G$-torsors on $K$ and the set $H^1(K,G)$, which sends the trivial 
$G$-torsor on $K$ to the trivial class in $H^1(K,G)$.
\end{thm}
Here, naturality means that the bijection is compatible with the changes 
of the field $K$ and the group $G$. For reader's convenience, we 
sketch a proof. 
\bp Let $\cG$ be a $G$-torsor on $K$. Fix a point $d_0\in\cG(\ov K)$. 
Then, for any $\s\in G_K$, there exists a unique $a_\s\in G(\ov K)$ 
such that $d_0^\s=d_0\, a_\s$ since $G_K$ acts on $G(\ov K)$ which 
itself is a group. It is not difficult to check that $\s\mapsto a_\s$ is an 
1-cocycle and hence induces an element in $H^1(K,G)$.

Conversely, if $\s\mapsto a_\s$ is an 1-cocycle, then on 
$G\times_K\ov K$, we obtain a new action of $G_K$ through 
$(\s,g)\mapsto \a_\s\,g^\s$. Since $G$ is quasi-projective, by Weil's 
theorem on descent, there exists a $K$-scheme $\cG$, or better, a 
$G$-torsor on $K$, such that $\cG\times_K\ov K=G\times_K\ov K$ is   
$G_K$-equivalent (with respect to the twisted action of $G_K$ on 
$G\times_K\ov K$). 
\ep

There is a canonical morphism $H^1(F,G)\to H^1(F_v,G)$ induced by the 
inclusion $F\hookrightarrow F_v$ for each $v\in S$, which themselves 
then induce a natural morphism
\be
H^1(F,G)\to \prod_{v\in S}H^1(F_v,G)
\ee
Denote by $\Ker^1(F,G)$ the kernel of this morphism. 
It is a result of Borel-Serre that $\Ker^1(F,G)$ is finite. 

\begin{cor} There are only finite many $G$-torsors 
on $F$, up to isomorphisms, such that, for all $v\in S$, the induced 
$G$-torsors on $F_v$  are trivial.
\end{cor}

\subsection{$G$-Torsors on $\ov{\Spec\,\O_F}$}

We begin with the following
\begin{defn} \normalfont
A {\it $G$-torsor on} $\ov{\Spec\,\O_F}$ is a scheme $\cG$, equipped with a flat 
surjective morphism $\pi:\cG\to \ov{\Spec\,\O_F}$ and a family of
flat surjective morphism $\pi_v:\cG_{\O_v^{~}}\to \ov{\Spec\,\O_v}$,
together with actions of $G(F_v)$ on $\cG_v$ for all 
$v\in S_\fin$, such that the induced morphism
\bes
\bm
\cG_x\times G_x&\lra& \cG_x\times_{\ov{\Spec\,\O_x}}\cG_x\\[0.7em]
(d,g)&\mapsto &(d,d^g)
\em
\ees
are isomorphisms for all points $x\in\Spec\, \O_F$, closed or 
generic.
\end{defn}

In particular, if $\cG$ is a $G$-torsor on $\ov{\Spec\,\O_F}$, the 
fiberwise $\{G(F_x)\}_{x\in X}$  acts on $\{\cG(F_x)\}_{x\in X}$ with 
respect to $\pi$, and the action on each generic fiber is faithful and 
transitive. Moreover, for each $v\in S_\fin$,  over the local integral base, 
induced by the natural inclusion $\O_F\hookrightarrow \O_v$, we 
obtain a composition of morphisms $\Spec\,\O_v\to \Spec\,\O_F\to 
\ov{\Spec\,\O_F}$. Thus, for the $G$-torsors $\cG_{\O_v}$ on 
$\Spec\,\O_v$, it is not too difficulty to deduce the following result, 
whose proof we left to the reader.

\begin{lem}\label{lem5}
\begin{enumerate}
\item [(1)] For each finite place $v\in S_\fin$,
induced from the natural morphisms $\O_v\hookrightarrow F_v$ and 
$\O_v\to k_v$, we have 
\bes
(\cG_{\O_v})_{\eta_v}\simeq  \cG_{\eta_v}\qqan (\cG_{\O_v})_{k_v}
\simeq  \cG_{k_v}.
\ees
\item[(2)] For each infinite place $\s\in S_\infty$, induced by the 
natural embedding $F\hookrightarrow F_\s$, we have
\bes
(\cG_\eta)_\s\simeq \cG_\s.
\ees
In particular, $(\cG_\eta)_\infty\simeq\cG_\infty$.
\end{enumerate}
\end{lem}

\subsection{Inner Form}
When working over integral base $\Spec\,\O_F$, our choice of a 
Chevalley basis $\big\{x_\a:\a\in\Phi\big\}\bigsqcup
\big\{h_\a:\a\in\De\big\}$ determines a pinning 
$(T,B,\{x_\a\}_{\a\in \De})$ of $G$. To deal with the associated 
compatibility problem, in the automorphism group $\Aut(G)$ of $G$, 
we consider the so-called {\it outer automorphism group} $\Out(G)$ 
defined as the collection of the automorphisms of $G$ which preserves 
the pinning $(T,B,\{x_\a\}_{\a\in \De})$. There is a natural split short 
exact sequence
\be\label{eq3}
1\to G^\ad\to\Aut(G) \buildrel{\pi}\over \to \Out(G)\to 1.
\ee
Indeed,  $G^\ad$ is identified with the image of $G$ under the adjoint 
representation of $G$,\footnote{The adjoint group $G^\ad$ of $G$ is 
the Zariski connected component of $\Aut(G)$, and may also be 
identified with the connected component of the automorphism group of 
$\frg:=\Lie(G)$. In particular, $G^\ad(F)$ coincides with  the group of 
inner automorphisms of G(F) defined over F, and acts simply transitively 
on the triples $(T,B,\{x_\a\})$ of all pinnings of $G(F)$.} and hence also 
fits into the short exact sequence
\be\label{eq4}
1\to G^\ad\to\Aut(G) \buildrel{\pi}\over \to \Aut(\Phi,\De)\to 1.
\ee
where $\Aut(\Phi,\De)$ denotes the automorphic of the root system 
$(\Phi,\De)$. In addition, the pinning $\{x_\a\}_{\a\in\De}$ identifies 
$\Aut(\Phi,\De)$ with $\Aut(G,T,B,\{x_\a\})$ and hence introduces a 
section $s: \Aut(\Phi,\De)\to \Aut(G)$ of the morphism 
$\pi: \Aut(G) \to \Aut(\Phi,\De)$ in \eqref{eq4}.

By taking Galois cohomology, we obtain the morphisms
\bes
H^1(F,\Aut(G))\buildrel{H^1_\pi}\over \lra H^1(F,\Out(G))\qan 
H^1(F,\Out(G))\buildrel{H^1_s}\over \lra H^1(F,\Aut(G)).
\ees
Hence, naturally associated to an element $\xi \in  H^1(F,\Aut(G))$ are 
the elements $H^1_\pi(\xi)\in H^1(F,\Out(G))$ and 
$(H^1_s\circ H^1_\pi)(\xi)\in H^1(F,\Aut(G))$. It is not too difficult to 
see that this element belongs to $H^1(F,G^\ad)$. Denote the induced 
$G^\ad$-torsor on $F$ by $G^\xi$, which, for later use, we call  {\it the 
inner form of} $G$ associated to $\xi$.

For example, an element $H^1(F,\Aut(G))$ induces naturally a $G$-
torsor on $F$. Hence, if we take $\xi$ as the trivial $G$-torsor, namely, 
$G$ itself to start with, then the element corresponding to 
$(H^1_s\circ H^1_\pi)(\xi)$ defines a split group $G^\ad$ on $F$.

\section{Compatible Metrics}\label{ComM}
In this subsection, we assume that our base field is the field of real numbers, unless otherwise stated explicitly.

\subsection{Maximal Compact Subgroup}\label{081003}

We first recall some basic facts on maximal compact subgroups of a 
real reductive group, mainly following  \cite{BS}.

Let $G$ be a real Lie group with finitely many connected components. 
Then any compact subgroup of $G$ is contained in a maximal compact 
subgroup. Moreover, if $K$ is a maximal subgroup of $G$, then $G$ is 
diffeomorphic to the direct product of $K$ with a euclidean space,  any  
maximal compact subgroup is conjugate to $K$ in $G$ and 
$G/G^0\simeq K/K^0$, where $\bullet^0$ denotes the connected 
component of $\bullet$ containing the unit element.

In addition, if $G_1$ is a closed normal subgroup of $G$ admitting only 
finitely many connected components, then the maximal compact 
subgroups of $G_1$ are the intersections of $G_1$ with maximal 
compact subgroups of $G$. Similarly, if $G_1$ is a closed subgroup 
of $G$ with finitely many  connected components such that all maximal 
compact subgroups of $G$ are conjugate by elements of $G_1$, the 
maximal compact subgroups of $G_1$ are the intersections of $G_1$ 
with the maximal compact subgroups of $G$. Consequently, in both 
cases, by taking a maximal compact subgroup $K$ of $G$ containing 
a maximal compact subgroup of $G_1$, we conclude that $G_1\cap K$ 
is a maximal compact subgroup in $G_1$ for at one and hence for all 
maximal compact subgroups of $G$ by conjugacy.

More generally, if $G\to G'$ is a surjective morphism (of Lie groups) 
whose kernel admits only finitely many connected components, then 
the maximal compact subgroups of $G'$ are the images of the maximal 
compact subgroups of $G$.

\subsection{The Cartan Involution}\label{thecart}
We here recall some basic facts on {\it the} Cartan involution 
associated to an algebraic group,  mainly following \cite{GV}.

Let $G$ be an algebraic group defined over a base field 
$F\subseteq \R$. Denote by $RG$ the {\it radical of $G$} and $R_uG$ 
the {\it unipotent radical  of $G$} and $R_dG$ the so-called {\it split 
radical of $G$}, namely, the greatest connected $k$-split subgroup of 
$RG$. By definition, a {\it Levi subgroup}\index{reductive group!
parabolic subgroup!Levi subgroup} of $G$ is a maximal reductive 
$k$-subgroup of $G$. Let 
\be
G^1=\bigcap_{\chi\in X(G)_k}\Ker(\chi^2),
\ee
with $X(G)_F:=\Hom_F(G,\G_m)$, the group of $F$-morphism of $G$ 
into $\G_m$.  Then $G^1$ is a normal subgroup of $G$, and is defined 
over $F$. Note that for a character $\chi$ in $X(G)_F$, its restriction to 
$G^1$ is of order $\leq 2$, hence is trivial on $(G^1)^0$. Consequently, 
\be
(G^1)^0=\big(\bigcap_{\chi\in X(G)_F}\Ker(\chi)\big)^0.
\ee 
Since any character in $X(G)_F$ is trivial on $R_uG$, we have 
$G^1=L^1\ltimes R_uG$ for any Levi subgroup $L$ of $G$. Hence, if 
$A$ is a maximal $F$-split torus of $RG$, then 
\begin{enumerate}
\item [(i)] $G(\R)=A(\R)^0\ltimes G(\R)^1$ and 
\item[(ii)] $G(\R)^1$ contains all compact subgroups of $G(\R)$. More 
generally,
\item[(iii)]  if $A_1$ and $A_2$ are two $F$-tori in $RG$ such that 
$A_1$ is $F$-split, $A_1A_2$ is a torus and $A_1\cap A_2$ is finite, 
then here exists a normal $F$-subgroup $G_1$ of $G$ containing 
$A_2$ and $G^1$ such that $G(\R)=A_1(\R)^0\ltimes G_1(\R)$.
\end{enumerate}

Therefore, if $P$ is a parabolic $F$-subgroup of $G$ and $A$ is a 
maximal $F$-split torus of the split radical $R_dG$ of $G$, then, for a 
maximal compact subgroup $K$ of $G(\R)$, we have 
\begin{enumerate}
\item [(a)] $K\cap P$ is a maximal compact subgroup of $P(\R)$, and 
\item[(b)] $G(\R)=KP(\R)=KA(\R)^0P(\R)$.  Furthermore, 
\item[(c)] if $KaP(\R)^0=Ka'P(\R)^0$ for some $a,\,a'\in A(\R)^0$, then 
$a=a'$ and the map $G(\R)\to A(\R)^0$ sending $g$ to $a=a(g)$ 
characterized by $g\in KaP(\R)^0$ is real analytic. 
\end{enumerate}

As a direct consequence, when $G$ is a reductive group, there exists 
one and only one involutive automorphism $\theta_K$ of $G(\R)$ 
associated to $K$ satisfying the following properties.
\begin{enumerate}
\item [(1)]  $\theta_K$ is "algebraic," i.e. the restriction to $G(\R)$ of an 
involutive automorphism of algebraic groups of the Zariski-closure of 
$G (\R)$ in $G$.\\[-0.88em]
\item[(2)] The fixed point set of $\theta_K$ is $K$.\\[-0.88em]
\item [(3)] If $G_1$ is a normal $\R$-subgroup of G, then 
$\theta_K(G_1(\R))=G_1(\R)$.\\[-0.88em]
\item [(4)]  $\theta_K$ leaves $[\frak g,\frak g]$ and $\frz$ stable.
Here, we use the same $\theta_K$ to denote the induced involution on 
$\frak g:=\Lie(G)$, and set $\frz$ denotes the center of $\frg$.
\\[-0.88em]
\item [(5)] If $\frak o$ is the (-1)-eigenspace of $\theta_K$ in $\frz$, 
then $V=\exp\frak o$ is a split component of $G$.\\[-0.88em]
\item [(6)] If $\frak p$ is the (-1)-eigenspace of $\theta_K$  in 
$\Lie(G(\R))$, then, for $\frak k:=\Lie(K)$, there is a decomposition 
$\Lie(G(\R))=\frak k\oplus \frak p$ and 
\bes
[\frak k,\frak k]\subseteq \frak k,\ \ [\frak k,\frak p]\subseteq \frak p,\ 
\ [\frak p,\frak p]\subseteq \frak k,\ \ \frak p^k\subset \frak p\ \ 
(\forall k\in K)
\ees
\item [(7)] The map $(k, X)\mapsto k\cdot \exp(X)$ is an isomorphism 
of analytic manifolds of $K\times\frak p$ onto $G(\R)$. 
\end{enumerate}

Note that in the case when $G$ is semi-simple, $\theta_K$ is the usual 
Cartan involution. Motivated by this, we call $\theta_K$ {\it the Cartan 
involution of $G(\R)$ with respect to $K$}\index{Cartan involution }.
Moreover, the existence of the Cartan involution $\theta_K$ implies an 
existence of a non-degenerate $G_\R$-symmetric bilinear form 
$\lan\cdot,\cdot\ran$ on $\frak g_\C\times \frak g_\C$ satisfying the 
follows.
\begin{enumerate}
\item [(a)] $\lan\cdot,\cdot\ran$ is invariant under $G$ and $\theta_K$, 
and is real on $\frak g\times\frak g$.\\[-1.40em]
\item [(b)] The quadratic form of $\lan\cdot,\cdot\ran$ is positive 
definite on $\frak p$ and negative definite on $\frak k$. In particular, if 
we set
\be\label{metc}
(X, Y) := -\lan X,\theta_KY\ran \qqan \|X\|^2 := -\lan X,\theta_KX\ran 
\quad\forall X, Y \in \frak g
\ee
then $(\cdot,\cdot )$ is a positive definite $K$-invariant and 
$\theta_K$-invariant scalar product for $\frak g\times\frak g$ with 
$\|\cdot\|$ its associated norm.\\[-1.40em]
\item [(c)] For $\frak g^1:=\Lie (G^1)$ and  $\lan\cdot,\cdot\ran^1$ 
and $\theta_K^1$ the restriction of $\lan\cdot,\cdot\ran$ and 
$\theta_K$ to $\frak g^1\times\frak g^1$ and $G^1$, respectively, we 
have that $(G^1,K,\theta_K^1,\lan\cdot,\cdot\ran^1)$ inherit all the 
properties of $(G,K,\theta_K,\lan\cdot,\cdot\ran)$ above.
\end{enumerate}
In addition,  since $\lan\cdot,\cdot\ran$ is $G$-invariant, the following 
infinitesimal invariance holds. 
\begin{enumerate}
\item [(d)] $\lan\cdot,\cdot\ran$ is characteristic, namely,  
\be\label{charB}
\ba
\lan \phi(X),\phi(Y)\ran=&\lan X,Y\ran\qquad \forall \phi\in\Aut_{\Lie}
(\frak g),\, \ \forall X,\, Y \in\frak g\\
\lan [X, Y],Z\ran =& \lan X,[Y,Z]\ran \qquad \forall X,\, Y,\,Z \in\frak g.
\ea \ee
\end{enumerate}

Therefore, $[\frak g,\frak g]$ and  $\frz$ are mutually orthogonal with 
respect to $\lan\cdot,\cdot\ran$. In addition, since 
$\lan\cdot,\cdot\ran$ is $\theta_K$-invariant, $\frak k$ and 
$\frak p$ is mutually orthogonal. Conversely, we may reconstruct  the 
bilinear form $\lan\cdot,\cdot\ran$ using all the above conditions. To 
be more precise, starting with the Cartan-Killing form on 
$[\frak g,\frak g]$, we may extend it to obtain $\lan\cdot,\cdot\ran$ as 
the direct sum of the Cartan-Killing form with a symmetric 
non-degenerate bilinear form on $\frz$, which is negative definite on 
$\frz\cap \frak k$ and positive definite on $\frz \cap \frak p$. Finally, 
we may  extend this latest $\lan\cdot,\cdot\ran$ to the total space 
$\frak g_\C\times \frak g_\C$.  For later use, we call 
$\lan\cdot,\cdot\ran=:\lan\cdot,\cdot\ran_K$ the {\it canonical form on 
$\frak g$ associated to $K$}\index{canonical form}. For later use, when 
$\lan\cdot,\cdot\ran_K$ is viewed as a linear form from $\frg$ to 
$\frg^*$, we write it as $H_K$.

Finally, let us point out that the Cartan involution can be applied in 
many ways. For example, if $G_1$ is a $\R$-subgroup of $G$ 
containing $R_uG$ such that all maximal compact subgroups $K$ of 
$G(\R)$ are conjugate under $G_1(\R)$, then, for a Levi subgroup $M$ 
of $G(\R)$, it makes sense to take about the Cartan involution 
$\theta_K$ of $M$ with respect to $K$. Moreover,  in this case, the 
subgroup $(G_1\cap M)\cap \theta_K(G_1\cap K)$ is the unique 
$\theta_K$-stable Levi subgroup of $G_1(\R)$ contained in $M$.  
Consequently, if $P$ is  a parabolic $\R$-subgroup of $G$, and $K$ is 
a maximal compact subgroup of $G(\R)$ and $M$ is a Levi subgroup 
of $G(\R)$ containing $K$, then $M\cap P$ contains one and only one 
tLevi subgroup of $P(\R)$ stable under $\theta_K$. For this reasons, 
we will fix a maximal compact subgroup $K$ of $G$ in the sequel.

\subsection{Fine Involutions for Maximal Compact Subgroups}
We are now ready to introduce new structures called fine involutions 
and their associated compatible metrics for general reductive groups, 
which may be viewed as natural generalizations of the known structures 
for semi-simple groups (see e.g. \cite{Gr2}). 

Let $G$ be a reductive group defined over a subfield $F\subset \R$. 
Fix a maximal compact subgroup $K$ of $G$. Motivated by the Cartan 
involution associated to the maximal compact group $K$ of $G$, we 
give the following:

\begin{lem}\label{le7} Let $H$ be a positive definite real symmetric 
bilinear form on $\frak g$. Assume $H$ is $K$-compatible with respect 
to the Lie structure of $\frak g$. Then, for $\theta_H:=-H_K^{-1}H$,
\begin{enumerate}
\item [(1)]  $\theta^2_H=1$.
\item [(2)] $H_KH^{-1}H_K=H$. That is to say, $H_K\theta_H=-H$.
\item[(3)] Let $\frak k_H$, resp. $\frak p_H$, be the $(+1)$-
eigenspace, resp. the $(-1)$-eigenspace, of $\theta_H$ on $\frak g$. 
Then
\begin{enumerate}
\item [(a)] $\frak g=\frak k_H\oplus \frak p_H$. 
\item [(b)] On $\frak k_H$, resp. $\frak p_H$, $H_K=-H$, resp. 
$H_K=H$, is negative definite, 
resp. positive definite. 
\item[(c)] $\frg=\frak k_H\oplus \frak p_H$ is an orthogonal 
decomposition with respect to $H_K$.
\end{enumerate}
\item [(4)] $H$ is compatible with $H_K$, i.e. $H_K:\frak g\to\frak g^*$ 
is an isometry with respect to the metric $H$ on $\frak g$ and  the 
metric $H^{-1}$ on $\frak g^*$.
\end{enumerate}
\end{lem}

\bp  (1) By the infinitesimal invariance of $H_K$, namely, two relations 
in \eqref{charB},  we have 
$\!~^t\theta_H H_K\theta_H=H_K$ since $\theta_H$ is compatible with 
the Lie structure on $\frg$. On the other hand, 
$\!~^t\theta_H H_K=-H_KH^{-1}H_K=H_K\theta_H$. Therefore 
$\theta^2_H=1$. 

(2) This is a direct consequence of (1). Indeed, since  $\theta^2_H=1$, 
we have $(H^{-1}H_K)(H^{-1}H_K)=\mathrm{Id}_\frg$. Therefore,  
$H_KH^{-1}H_K=H$ and hence $H_K^theta_H=-H$.

(3) (a) This is a standand result in linear algebra.

(b)   This is a direct consequence of (2). Indeed, since 
$H_K\theta_H=-H_KH^{-1}H_K=-H$, we have that $H_K=-H$ on $\frk$, 
resp. $H_K=H$ on $\frp$, by the fact that, on $\frk$, resp. on 
$\frp$, $\theta_H=1$.  Consequently, $H_K$ is negative definite on 
$\frk$ and positive definite on $\frp$. 

(c) This is a direct consequence of (b) and (c). 

(4) This is a reinterpretation of (2) and (3). Indeed, by (3), 
$H_K:\frak g\to\frak g^*$ is an isomorphism. Moreover, by (2), 
$H_KH^{-1}H_K=H$, we obtain the following commutative diagram of 
isomorphisms
\bes
\bm 
\quad\frg&\buildrel{H_K}\over\lra&\frg^*\quad\\[0.3em]
~^H\!\xdownarrow&&\quad\xdownarrow \!\!~^{H^{-1}}\\
\quad\frg^*&\buildrel{H_K}\over\lla&\frg.\quad\\
\em
\ees
This implies that  $H_K:\frak g\to\frak g^*$ is an isometry with respect 
to the metric $H$ on $\frg$ and the metric $H^{-1}$ on $\frak g^*$.
\ep

\begin{defn}  \normalfont
An element $\theta$ in the automorphisms group $\Aut_\Lie(\frg)$ of 
the Lie algebra $\frg$  is called a {\it fine involutions} of $K$ (with 
respect to the Lie structure on $\frg$) if there exists a 
positive definite real symmetric bilinear form  $H$ on $\frg$ such that 
$\theta=-H^{-1}H_K$ and satisfies all the properties (1), (2), (3) and (4) 
in Lemma\,\ref{le7}. Here,  $H_K:\frg\to\frg^*$ denotes the linear 
isomorphism associated to the bilinear form $\lan\cdot,\cdot\ran_K$ 
on $\frg$ induced from the Cartan involution $\theta_K$ associated 
with $K$. Moreover, if this is the case, $H$ is called a 
{\it $K$-compatible metric on $\frg$ (with respect to its Lie structure)},
\footnote{Here as usual, we view the bilinear $H$ as a linear map from 
$\frak g$ to $\frak g^*$. Hence $H^{-1}H_K$ is indeed a linear 
endomorphism of $\frak g$.} and we denote $\theta$ by $\theta_H$.  
\end{defn}

From the discussion above, it is not difficult to see that fine involutions 
and admissible metrics on $\frg$ associated to $K$ works exactly in 
the same way for  $G^1$ as well, since $G^1$ is reductive and all 
maximal compact subgroups of $G$ are contained in $G^1$. Indeed, 
the corresponding constructions on $G^1$ may be viewed as the 
restrictions of the structures from $\frg$ to $\frg^1:=\Lie(G^1)$. For 
later use, set
\be\label{eq9}
\frg=\frg^1\oplus \frv.
\ee

\subsection{Compatible Metrics for Maximal Compact Subgroups}
\label{s1.2.4}

Denote by $\cM_{\frak g;K}^\tot$, resp. $\cM_{\frg^1;K}^\tot$ the 
moduli space of $K$-compatible euclidean metrics on $\frg$, resp. on 
$\frg^1$. Since they contains (the isometric class of) 
$(\cdot,\cdot)_K$ in \eqref{metc},  both $\cM_{\frak g;K}^\tot$ and 
$\cM_{\frg^1;K}^\tot$ are not empty. Moreover, $\cM_{\frak g;K}^\tot$, 
resp. $\cM_{\frg^1;K}^\tot$,  admits a natural interpretation as a 
subspace of the space of (isomety classes of) euclidean metrics on 
$\frak g$, resp. on $\frg^1$. Our main result of this section is the 
following:

\begin{prop}\label{prop8} Let $G$ be a reductive group defined over a 
subfield $F\subset \R$ and let $K$ be a maximal compact subgroup of 
$G$. Set $\frg=\Lie(G)$. We have
\begin{enumerate}
\item [(1)] There are natural actions of $G$ on $\cM_{\frak g;K}^\tot$ 
and $\cM_{\frg^1;K}^\tot$.
\item [(2)] The action of $G$ on $\cM_{\frak g;K}^\tot$ in (1) induces a 
natural diffeomorphism
\be
G^1(\R)/K\simeq \cM_{\frg^1;K}^\tot
\ee 
\end{enumerate}
\end{prop}
\bp (1) Recall that, for $\phi\in\Aut_\Lie(\frak g)$, we have 
$\!~^t\phi H_K\phi=H_K$ since $H_K$ is characteristic by 
\eqref{charB}. Hence, for any $H\in \cM_{\frak g;K}^\tot$, we have 
\bes
(-\!~^t\phi H\phi)H_K=\phi^{-1}\theta\phi.
\ees 
This implies that $\!~^t\phi H\phi$ belongs to 
$\cM_{\frak g;K}^\tot$ as well. Consequently,  the assignment 
\be
(H,g)\longmapsto \!~^t(\Ad\,g)H(\Ad\,g)\qquad
\forall H\in \cM_{\frak g;K}^\tot,\ \forall g\in G
\ee
defines a natural action of $G$ on $\cM_{\frak g;K}^\tot$. Here, as 
usual, for $g\in G$, $\Ad\,g$ denotes its adjoint. This proves (1).

(2) In terms of the fine involutions $\theta=\theta_H$, the action above 
is equivalent to 
\be
(\theta,g)\longmapsto \!~^t(\Ad\,g)\theta(\Ad\,g)\qquad\ \forall g\in G
\ee
Recall that, for a general $G$, if we set 
\bes
\frv_{\theta}:=\left\{\,X\in \frz: \theta X=-X\,\right\}.
\ees 
Then $V_\theta:=\exp(\frv_\theta)$ is a split component of $G$.
\footnote{Denote by $\Ad: G\to\Aut_\Lie (\frg)$ the adjoint action of 
$G$ on $\frg$. Then $\Ker(\Ad)$ acts trivially on $\frg$. A {\it split 
component} of $G$ is defined to be a maximal closed linear subspace 
of $\Ker(\Ad)$. By Proposition 2.1.5 of \cite{GV}, if $V$ is a split 
component of $G$, then $G=G^1\cdot V$.} Moreover, by Proposition 
2.1.10 of \cite{GV}, the assignment $\theta\mapsto \frv_\theta$ gives a 
bijection from the set of fine involutions associated to $K$ to the set of 
split components of $G$. In particular, when $G=G^1$, this  map gives a 
bijection from the set of fine involutions to the set of maximal 
compact subgroups of $G$. Hence, in our case, since $G$ is connected 
and reductive, all maximal compact subgroups of $G$ are conjugate to 
each other. This implies that $G$ and hence $G^1$ acts transitively on 
$\cM_{\frak g^1;K}^\tot$. Thus, to complete our proof, it suffices to 
show that the stabilizer group of $\theta_K$ in $G^1$ is exactly $K$ 
itself.  For this we choose $\Theta_K:G^1(\R)\to G^1(\R)$ to be a 
Cartan involution satisfying $d\Theta_K=\theta_K$. By 
definition, 
\bes
K=\left\{\,g\in G^1(\R): \Theta(g)=g\,\right\}.
\ees 
On the other hand, for $g\in G^1$, 
$\theta_K=(\Ad g)^{-1}\theta_K(\Ad g)$ if and only if 
$g^{-1}\Theta_K g$ is in the center of $G(\R)^0$, the connected 
component of $G(\R)$ containing the unit element. But this center is 
trivial by our assumption, hence $g$ belongs to the stabilizer 
group of $\theta_K$ if and only if $g\in K$.
\ep

From the proof, we conclude that $\frv$ in \eqref{eq9} is identified with 
$\frv_H$ for a certain compatible $H$ of $K$.  As a direct 
consequence, we obtain the following:
\begin{cor}\label{cor9} Denote by $\cM_\frv$ be the moduli space of 
euclidean metric on $\frv$. Then
\begin{enumerate}
\item [(1)] $\cM_\frv\simeq\GL_{\dim_\R\frv}(\R)/O_{\dim_\R\frv}(\R)$, 
where $O_n(\R)$ denotes the orthogonal group of degree $n$.
\item[(2)] The natural map defined by
\bes
\bm
\cM_{\frak g;K}^\tot&\lra&\cM_{\frak g^1;K}^\tot\times \cM_\frv
\\[0.6em]
H&\mapsto&(H|_{\frg^1}, H|_\frv)
\em
\ees
is bijecive.
\end{enumerate}
\end{cor}

\section{Arithmetic $G$-Torsors on $\ov{\Spec\,\O_F}$}
\subsection{Integral Structures on Lie Algebras}

Let $F$ be an algebraic number field with $\O_F$ the ring of integers, 
and let $G$ be a connected  split reductive group over $F$ with $\frg$ 
its Lie algebra. Our aim here is to introduce a $G(\O_F)$-invariant 
projective $\O_F$-module $\frak g^{}_{\O_F}$ in $\frak g\otimes_F\R$ 
which is closed under the Lie operation.

For simplicity, assume $F=\Q$. Since $G$ is defined over $\Q$, its Lie 
algebra $\frak g$ admits a natural rational structure $\frak g_\Q$ and 
the adjoint representation $G\to\End_\Lie(\frak g_\Q)$ is a morphism 
defined over $\Q$. Consequently, there always exist $G(\Z)$-invariant 
integral structures in $\frak g_\Q$, since, for any integral structure in 
$\frak g_\Q$, the image under the action of $G(\Z)$ is again an 
integral structure in $\frak g_\Q$. Obviously, the summation of two 
$G(\Z)$-invariant integral structures in $\frak g_Q$ is again a 
$G(\Z)$-invariant integral structure. Moreover, if $\frak g_\Z$ is a 
$G(\Z)$-invariant integral structure in $\frak g_Q$, we have 
$[\frak g_\Z,\frak g_Z]\subset\frak g_\Q$. Hence, by clearing up 
denominators, we  can instead assume that  
$[\frak g_\Z,\frak g_Z]\subseteq\frak g_\Z$ from the beginning. In this 
way, we obtain a unique maximal $G(\Z)$-invariant integral structure 
$\frak g_\Z$ in $\frak g_Q$ satisfying the condition that 
$[\frak g_\Z,\frak g_Z]\subseteq\frak g_\Z$. 

Put this in a more concrete form, since our reductive group $G$ is 
defined over $\Q$, we may use the structural decomposition 
\bes
\frg_\Q=\frg^{\ses}_\Q\oplus \frz_\Q
\ees
where $\frg^{\ses}_\Q$ denotes the rational structure on $\frg^\ses$ 
induced by the semi-simple Lie sub-algebra of $\frg$. Since 
$\frv\subset \frz$, this decomposition is compatible 
with $\frg_\Q=\frg^1_\Q\oplus \frv_\Q$ induced by \eqref{eq9}. 
Moreover, since $\frz$ is abelian,  we obtain a natural decomposition
\be\label{eq13}
\frz_\Q=\frv_\Q\oplus\frg^1_\Q/\frg_\Q^\ses.
\ee
Now  by applying the Chevalley (integral) basis for semi-simple Lie 
algebras, we obtain a canonical integral structure  $\frg^\ses_\Z$ on 
$\frg_\Q^\ses$. Hence, what is left is  to introduce an integral structure 
on $\frz_\Q$ which is compatible with the decomposition 
\eqref{eq13}. But this is trivial since $\frz$ is an abelian sub Lie algebra 
defined over $\Q$. We thus obtain an induced integral structure on 
$\frz_\Q$, which we denoted  by $\frz_\Z$. Similar arguments then lead 
to the integral structures $\frg^1_\Z$ and $\frv_\Z$ on $\frg^1_\Q$ 
and $\frv_\Q$, respectively. Consequently, we have 
\be
\frg_\Z:=\frg_\Z^\ses\oplus \frz^{~}_\Z=\frg_\Z^1\oplus \frv^{~}_\Z.
\ee
 
The above discussion works well if we replace $G/\Q$ by a split 
reductive group $G/F$ with $F$ a general number field. To indicate the 
dependence on $F$, we rewrite the associated Lie algebra by $\frg_F$. 
Since it admits a natural $F$-linear space structure, through the 
Minkowski embedding 
$F\hookrightarrow F_\infty:=\prod_{\s\in S_\infty}F_\s$, we obtain a Lie 
algebra
\be
\frak g_\infty:=
\frak g_F\otimes_\Q\R:=\prod_{\s\in S_\infty}\frak g\otimes F_\s.
\ee 
We introduce an $\O_F$-lattice structure on $\frg_F$ by setting
\be
\frak g^{}_{\O_F}:=\frg_\Z\otimes\O_F\hookrightarrow 
\frg_F\hookrightarrow \frak g_\infty.
\ee
Then
\be
\frak g^{}_{\O_F}=\frg_{\O_F}^\ses\oplus \frz^{~}_{\O_F}
=\frg_{\O_F}^1\oplus \frv^{~}_{\O_F}
\ee
where $\frg_{\O_F}^\ses:=\frg_\Z^\ses\otimes\O_F,\ \frz^{~}_{\O_F}
=\frz^{~}_\Z\otimes\O_F,\ 
\frg_{\O_F}^1=\frg_\Z^1\otimes\O_F$ and $\frv^{~}_{\O_F}
= \frv^{~}_\Z\otimes\O_F$.
Using a similar argument as above, we conclude that 
$\frak g^{}_{\O_F}$ is a projective 
$\O_F$-submodule in $\frak g_F$ such that 
$[\frak g^{}_{\O_F},\frak g^{}_{\O_F}]\subseteq\frak g^{}_{\O_F}$.  
In the sequel, $\frak g^{}_{\O_F}$ will be called {\it the canonical 
infinitesimal $\O_F$-structure} of $G/F$.

\subsection{Arithmetic $G$-Torsors}

Let $G$ be a split reductive group over a number field $F$. For each 
$\s\in S_\infty$, we fix a 
maximal compact subgroup $K_\s$ of $G^1(F_\s)$, set 
$r_{G,\s}^{G^\ses}:=\rk(G)-\rk(G^\ses)$.  We denote a real, resp. 
complex, $\s\in S_\infty$ by $\s:\R$, resp. $\s:\C$. By \S\ref{s1.2.4}, 
we obtain the moduli spaces 
$\cM_{\frak g^{}_{F_\s}\!;K_\s}^\tot$, resp. 
$\cM_{\frak g^1_{F_\s}\!;K_\s}^\tot$, of the compatible metrics with 
respect to $K_\s$ on $\frak g^{}_{F_\s}$, resp. $\frak g^1_{F_\s}$, and 
natural isomorphisms 
\bea
\prod_{\s\in S_\infty}\cM_{\frak g^1_{F_\s}\!;K_\s}^\tot\stackrel{\simeq}
{\longrightarrow}& 
\prod_{\s\in S_\infty}G^1(F_\s)/K_\s,\\
\prod_{\s\in S_\infty}
\cM_{\frak g^{}_{F_\s}\!;K_\s}^\tot\stackrel{\simeq}{\longrightarrow}& 
\prod_{\s\in S_\infty}\cM_{\frak g^1_{F_\s}\!;K_\s}^\tot\\
\times&\Big(\!\big(  \GL_{r_{G,\s}^{G^\ses}}(\R)/O_{r_{G,\s}^{G^\ses}}
(\R))^{r_1}\!\!\times\!
\big( \GL_{r_{G,\s}^{G^\ses}}(\C)/U_{r_{G,\s}^{G^\ses}}(\C)\big)^{r_2}\!
\Big).
\eea
Here $O_n(\R)$, resp. $U_{n}(\C)$, denotes the orthogonal group, resp. 
the unitary group, of degree $n$, and $r_1$, resp. $r_2$, denotes the 
number of real, resp. complex, places of $F$. For later use, set 
\be
G^1(F_\infty)/K(F_\infty):=\!\!\prod_{\s\in S_\infty}\!\! 
G^1(F_\s)/K_\s\qan \cM_{\frak g^\bullet_\infty;K_\infty}^\tot\!\!
:=\!\!\prod_{\s\in S_\infty}\!\!\cM_{\frak g^{\bullet}_{F_\s}\!;K_\s}^\tot
\ee
where, to simplify our notations, we use $\frg^\bullet$ as a running 
symbol for $\frg$ and $\frg^1$.

\begin{defn} \normalfont
\normalfont  Let $G/F$ be a connected (split) reductive group and let 
$K_\infty:=(K_\s)_{\s\in S_\infty}$ be a family of maximal compact 
subgroups of $(G(F_\s))_{\s\in S_\infty}$. By a 
{\it $K_\infty$-compatible arithmetic $G$-torsor} over 
$\ov{\Spec\,\O_F}$, or simply over $\O_F$, we mean a tuple 
$(\cG, (H_\s)_{\s\in S_{\infty}})$ consisting of a $G$-torsor $\cG$ on 
$\ov{\Spec\,\O_F}$ and  an element $(H_\s)_{\s\in S_\infty}$ of 
$\cM_{\frak g_\infty;K_\infty}^\tot$. 
\end{defn}

Even apparently not quite related, $H_\infty:=(H_\s)_{\s\in S_\infty}$ may be 
viewed as a family of $(K_\s)_{\s\in S_\infty}$-compatible 
metrics on the tangent bundles of the $G$-torsor 
$\cG_\infty:=\prod_{\s\in S_\infty}\cG_\s$. To explain this, we first 
recall that $\cG_\eta$ admits a natural $G(F)$-torsor structure. This, via 
the Minkowski embedding, induces a natural $G(F_\infty)$-torsor 
structure on $\cG_\infty$. Hence, for a fixed base point of 
$\cG_\infty$, the the tangent space of $\cG_\infty$ at this point is 
canonically identified with $\frg_\infty$. Consequently, we obtain a 
natural metric on this tangent space. Moreover, since  $\cG_\infty$ is 
a $G(F_\infty)$-torsor, its tangent bundle is a flat bundle. Thus, with the 
help of the so-called parallel transforms, we obtain a natural metric on 
the tangent bundle of $\cG_\infty$. This metric is uniquely determined 
by $H_\infty$.

Furthermore, working over $G^1$, at infinite places, $H_\infty$ is a 
$K_\infty$-compatible metric on $\frg^1_\infty$. Thus 
$(\frg^1_{\O_F},H_\infty)$ for an $\O_F$-lattices in $\frg^1_\infty$ 
whose projective $\O_F$-module component  $\frg^1_{\O_F}$ 
also admits a natural Lie structure over $\O_F$ because, by our 
construction,
$[\frg^1_{\O_F},\frg^1_{\O_F}]\subseteq \frg^1_{\O_F}$.

\begin{defn} \normalfont
Let $(\La,(H_\s))$ be a pair consisting of a projective $\O_F$-module 
$\La\subset \frg_F$ and a family of $K_\infty$-compatible metrics on 
$\frg_\infty$. If $\La$ is $G(\O_F)$-invariant under the adjoint action, 
and $[\La,\La]\subset\La$, then $(\La,(H_\s))$ is called a 
{\it $K_\infty$-compatible principle $G^1$-lattices over $\O_F$}.
\end{defn}

Denote by $\cM_{G,F}^\tot$, resp. $\cM_{G^1,F}^\tot$, the moduli stack 
of $K_\infty$-compatible arithmetic $G$-torsors, resp. $G^1$-torsors, 
on $\ov{\Spec\,\O_F}$. Then we have the following

\begin{thm} Let $G$ be a split reductive group on $F$. Then there exist 
the following natural identifications:
\begin{enumerate}
\item [(1)] $\cM_{G^1,F}^\tot\simeq \prod_{\xi\in\Ker^1(F,G^1)}
G^{1,\xi}(F)\backslash G^1(\A)/\prod_{v\in S_\fin}G^1(\O_v)\times 
K_\infty$.
\item [(2)] $\cM_{G,F}^\tot\simeq \dis{\prod_{\xi\in\Ker^1(F,G)}}\!\!\!
G^{\xi}\!(F)\backslash G(\A)/
\big(\!\!\!\!\prod_{v\in S_\fin}\!\!\!\!G(\O_v)\times K_\infty\!\times 
O_{r_{G,\s}^{G^\ses}}\!(\R)^{r_1}\!
\times U_{r_{G,\s}^{G^\ses}}\!(\C)^{r_2}\!\big).$
\end{enumerate}
\end{thm}
\bp
(1) For each $v\in S_\fin$, let $X_v=\Spec(\O_v)$ and set $X_v^\bullet$ 
be the complementary open subset of $\{v\}$ in $X_v$. Then for each 
element $g_v$ of the affine Grassmannian  $G(F_n)/G(\O_v)$, we obtain 
a $G(F_v)$-torsor $\cE_v$ on $X_v$ equipped with a trivialization 
on $X_v^\bullet$ with the trivial $G$ torsor, so that if an automorphicm 
of $\cE_v$ is trivial on $X_v^\bullet$, then it is necessarily trivial. Recall 
that  there exists a natural morphism 
$G(F_v)/G(\O_v)\lra G(F)\backslash \prod'G(F_v)/G(\O_v)$ which maps 
$g_v\in G(F_v)/G(\O_v)$ on the tuple consisting of $g_v$ at $v$ and 
the unit elements of $G(F_{v'})/G(\O_{v'})$ for all places 
$v'\in S_\fin\!\smm\!\{v\}$. Then, by the compatibility condition 
in Lemma\,\ref{lem5}, and the fact that $\Spec\,\O_F$ is affine and 
$\O_F$ is Dedekind domain, we conclude that the local maps induces 
an identifications of $G(F)\backslash\prod'G(F_v)/G(\O_v)$ with the 
moduli spaces of $G^1$-torsors on $\ov{\Spec\,\O_F}$, by adopting a 
result of Beauville-Laszlo \cite{BZ} in the case when $G=\GL_n$ and 
that of Heinloth in \cite{He} for general cases. Therefore, to conclude 
our proof, it suffices to apply Proposition\,\ref{prop8} to take care of 
the factor of $K_\infty$-compatible metrics $(H_\s)$, since $G(\A_F)$ 
is nothing but $\prod_{v\in S}'G(F_v)$. 
 
(2) is a direct consequence of (1) and its proof, if we apply Corollary\,
\ref{cor9}.
\ep

In the sequel, to simplify our presentations, we will use arithmetic 
$G$-torsors instead of the full version of $K_\infty$-compatible  
arithmetic $G$-torsors over $\ov{\Spec\,\O_F}$, if no confusion arises.

\subsection{Slopes of Arithmetic $G$-Torsors}

Let $T_G$ be the maximal split torus in the center $Z_G$ of $G$ and  
let $T_G'$ be the maximal split quotient torus of $G$. Then
\bes
T_G=\Hom(\G_m,Z_G)\otimes \G_m\qqan 
T_G'=\Hom(\Hom(G^\ab,\G_m),\G_m).
\ees
Here $G^\ab:=G/[G,G]$ denotes the maximal abelian quotient of $G$. It 
is well known that the composition 
\bes
T_G\hookrightarrow Z_G\hookrightarrow G
\twoheadrightarrow G^\ab\twoheadrightarrow T_G'
\ees
is an isogeny, i,.e. a morphism with finite kernel and cokernel. 
Consequently, if we set
\bea
X_*(T_G):=\Hom(\G_m,T_G),\qquad X_*(T_G'):=\Hom(\G_m,T_G')
\\[0.5em]
X^*(T_G):=\Hom(T_G,\G_m),\qquad X_*(T_G'):=\Hom(T_G',\G_m)
\eea
then $X_*(T_G)\hookrightarrow X_*(T_G')$ is a free abelian group of the 
same rank. Moreover, since $\Hom(\G_m,\G_m)\simeq \Z$, there is a 
non-degenerate pairing
\be\label{eq21}
\bm
\lan\cdot,\cdot\ran:&X_*(T_G^\bullet)\times X^*(T_G^\bullet)&\lra&\Z
\\[0.5em]
&(\chi,\mu)&\mapsto&\lan\chi,\mu\ran
\em
\ee
Here $T_G^\bullet= T_G$ or $T_G'$. Set now 
\be
X_*(T_G^\bullet)_\infty:=\prod_{\s\in S_\infty}X_*(T_G^\bullet)^{~}_{F_\s}
\qan
X^*(T_G^\bullet)_\infty:=\prod_{\s\in S_\infty}X^*(T_G^\bullet)^{~}_{F_\s},
\ee
then \eqref{eq21} induces a non-degenerating pairing
\be
\bm
\lan\cdot,\cdot\ran:&X_*(T_G^\bullet)_\infty\times 
X^*(T_G^\bullet)_\infty&\lra&\C\\[0.5em]
&(\chi,\mu)&\mapsto&\lan\chi,\mu\ran.
\em
\ee
Obviously, there is a natural action of 
$\prod_{\s\in S_\infty}\Gal(F_\s/\R)$ on $X^*(T_G^\bullet)_\infty$, and 
\be
\left\lan X_*(T_G^\bullet), X^*(T_G^\bullet)_\infty^{\prod_{\s\in S_\infty}
\Gal(F_\s/\R)}\right\ran\subseteq\R.
\ee
For our own convenience, we denote the invariant space 
$X^*(T_G^\bullet)_\infty^{\prod_{\s\in S_\infty}\Gal(F_\s/\R)}$ by 
$X^*(T_G^\bullet)_\infty^\ar$.

\begin{defn} \normalfont 
Let $\ov\cG=(\cG,(H_\s))$ be a $K_\infty$-compatible arithmetic $G$-
torsor over $\ov{\Spec\,\O_F}$. An element $\mu\in X_*(T_G')_\infty$ is 
called the  {\it slope of $\ov\cG$}, denoted 
by $\mu(\ov\cE)_\infty$, if, for all $\chi\in X^*(A_G')$, we have
\be
\lan\chi,\mu\ran=\deg_{\ar}(\ov \cG_\chi)
\ee
where $\ov \cE_\chi$ denotes an arithmetic $\G_m$-torsor on 
$\ov{\Spec\,\O_F}$ induced by the reduction of structure group 
$G\twoheadrightarrow T_G'\buildrel{\chi}\over\to \G_m$, and 
$\deg_{\ar}$ denotes the arithmetic degree.
\end{defn}

This definition makes sense, since an arithmetic $\G_m$-torsor on 
$\ov{\Spec\,\O_F}$ is simply a metrized line bundle on 
$\ov{\Spec\,\O_F}$. Hence  its arithmetic degree 
$\deg_\ar(\ov \cG_\chi)$ is well-defined.

\begin{rem}\normalfont
\begin{enumerate}
\item [(1)] Arithmetic $G$-torsors are first introduced in my book \cite{W}. There 
is a serious overlap between this section and \S 16.2 of \cite{W}, 
even the context here is much clearer.  
\item[(2)]
As to be expected, the slope can be use to definite stability of 
arithmetic $G$-torsors (\cite{W}). We omit the details, since it will not 
be used in our current work.
\end{enumerate}
\end{rem}

\section{Arithmetic Characteristic Curves}
Let $G/F$ be a split reductive group with pinning
$(T,G,\{x_\a\}_{\a\in\De})$ such that $\{x_\a\}_{\a\in\De}$ can be extended to a Chevalley basis of $\frg=\Lie(G)$. 

Let $(\cG,H)$ be an arithmetic $G$-torsor on $\ov X$. Denote by 
$\frg_X$ the induced locally free sheaf on $X$, or equiva;ently,
the associated projective $\O_F$-module in $\frg_\infty$.
 
Let $(\cL,\rho)$ be a metrized invertible sheaf on $\ov X$. Denote by
$L$ its associated rank one projective $\O_F$-module. Since $\O_F$ is 
Dedekind, we may identify $L$ with a certain fractional ideal of $F$.  
Denote this fractional ideal by the same letter $L$, by an abuse of 
notations.

For an element $\varphi\in \frg_X\otimes_{\O_F}L\subset\frg_F\subset 
\frg_\infty$, we denote its images under the Chevalley characteristic 
morphisms $\chi^{~}_F$ and $\chi^{~}_\infty$ by 
$\chi(\varphi_F)\in \frc_F$ and $\chi(\varphi_\infty)\in \frc_\infty$, 
respectively. It is not difficulty to see that $\chi(\varphi)\in \frc({\O_F})$. 

The element $\varphi$ can also be used to identify $X$ with the 
horizontal section associated to $\varphi\,\O_X$ in 
$\frg_X\otimes\cL$. In this way, we obtain a morphism 
$a_\varphi: X\to\frc_{\O_F}$.

\begin{defn} \normalfont
Let $G/F$ be a split reductive group.
\begin{enumerate}
\item [(1)] The pair $(\cG,\varphi)$ consisting of an arithmetic 
$G$-torsor $\cG$ and an element 
$\varphi\in \frg_X\otimes_{\O_F}\cL$ with $\cL/\Spec\,\O_F)$ an invertible line sheaf is called an arithmetic Higgs $G$-torsor. 
\item[(2)] The {\it characteristic arithmetic curve} associated to 
$(\cG;\cL;\varphi)$ is defined to be the scheme $X_\varphi$ of 
arithmetic dimension one obtained from the morphism 
$X\buildrel{a_\varphi}\over\to \frc_{\O_F}^{~}$ through the  base 
change $\frt_{\O_F}^{~}\to \frc_{\O_F}^{~}$. That is to say, 
\bes
X_\varphi=X\times_{\frc_{\O_F}^{~}}\frt_{\O_F}^{~}
\ees
induced from the product diagram
\be
\bm
X_\varphi&\buildrel{\chi_\varphi^{~}}\over\lra&X\\
\xdownarrow&&\quad\,  \xdownarrow\, a_\varphi\\
\frt_{\O_F}^{~}&\buildrel{\chi}\over\lra&\frc_{\O_F}^{~}.
\em
\ee
\end{enumerate}

\end{defn}

Obviously,  the morphism $\chi_\varphi:X_\varphi\to X$ is a finite 
$(|W|:1)$-covering, even highly ramified in general. Here, as usual, $W$ 
denotes the Weyl group of $G/F$.

\begin{rem} \normalfont
\begin{enumerate}
\item [(1)] The above construction is motivated by the construction
of spectral curve (for $G=\GL_n$) and cameral curve (for general 
reductive $G$) by Beauville-Narasimhan (\cite{BN}) and Donagi-
Gaistgory (\cite{DG}), respectively.  
\item[(2)] When $G=\GL_n$, with the identification 
$\frt \simeq\Spec\,F[\frt]$, the Chevalley characteristic morphism may 
be viewed as the assignment for diagonal matrices to their unorded 
eigenvalues. For this reason, we sometimes also call $X_\varphi$ the 
{\it arithmetic eigen curve} of $X$ associated to $\varphi$. 
\end{enumerate}
\end{rem}

In the forthcoming papers (\cite{W2}, \cite{W3}), we will use 
arithmetic characteristic curves to construct arithmetic Hitchin 
fibrations and study the intersection homologies and 
perverse sheaves for the associated structures, following (Laumon-)Ngo's approach to the 
fundamental lemma (\cite{Ngo}) using Hitchin fibrations (\cite{Ngo1}).

~\vskip 8.0cm

Lin WENG

Institute for Fundamental Research

The $L$\,Academy 

{\it and}

Faculty of Mathematics 

Kyushu University

Fukuoka, 819-0395

JAPAN

E-Mail: weng@math.kyushu-u.ac.jp

\begin{thebibliography}{80}


\bibitem{BN} A. Beauville, M.S. Narasimhan, S. Ramanan, Spectral 
curves and the generalised theta divisor. J. Reine Angew. Math. 398 
(1989), 169-179.

\bibitem{BS} A. Borel, J.P. Serre,  Corners and arithmetic groups.  
Comment. Math. Helv. 48 (1973), 436–491.

\bibitem{DG} R.Y. Donagi, D. Gaitsgory,  The gerbe of Higgs bundles. 
Transform. Groups 7 (2002), no. 2, 109-153.

\bibitem{Gr2} D.R. Grayson,  Reduction theory using semistability. II. 
Comment. Math. Helv. 61 (1986), no. 4, 661-676.

\bibitem{GV}
R. Gangolli, V.S. Varadarajan, {\it Harmonic analysis of spherical 
functions on real reductive groups}. Ergebnisse der Mathematik und 
ihrer Grenzgebiete, 101. Springer-Verlag, Berlin, 1988.

\bibitem{H} J.E. Humphreys, {\it Introduction to Lie algebras and 
representation theory}. Second printing, revised. Graduate Texts in 
Mathematics, 9. Springer-Verlag, New York-Berlin, 1978.

\bibitem{NSW} J. Neukirch, A. Schmidt, K. Wingberg,  {\it Cohomology of 
number fields}. Grundlehren der Mathematischen Wissenschaften 323. 
Springer-Verlag, Berlin, 2000.

\bibitem{Ngo1} Fibration de Hitchin et endoscopie. Invent. Math. 164 
(2006), no. 2, 399-453.

\bibitem{Ngo} B.C. Ngô, Le lemme fondamental pour les algebres de 
Lie.  Publ. Math. Inst. Hautes Etudes Sci. No. 111 (2010), 1-169. 

\bibitem{W} L. Weng, {\it Zeta Functions of Reductive Groups and Their 
Zeros},  World Scientific 2018.

\bibitem{W2} L. Weng, Arithmetic Hitchin Fibrations I: Basic Properties, in preparations

\bibitem{W3} L. Weng, Arithmetic Hitchin Fibrations II: Perverse Sheaves and Intersection Homology, in preparations
\end{thebibliography}
\end{document}